\makeatletter% do not load the package program <<< added <<<<<<<<<<<<<<<<
\disable@package@load{program}{}
\makeatother

\documentclass[pdflatex,sn-mathphys-num,iicol]{sn-jnl}% Math and Physical Sciences Numbered Reference Style 

\usepackage{graphicx}%
\usepackage{multirow}%
\usepackage{amsmath,amssymb,amsfonts}%
\usepackage{amsthm}%
\usepackage{mathrsfs}%
\usepackage[title]{appendix}%
\usepackage{xcolor}%
\usepackage{textcomp}%
\usepackage{manyfoot}%
\usepackage{booktabs}%
\usepackage{algorithm}%
\usepackage{algorithmicx}%
\usepackage{algpseudocode}%
\usepackage{listings}%
\usepackage{color,soul}
\usepackage{lmodern}

\usepackage{tikz}
\usepackage{orcidlink}

\jyear{2021}%

\usepackage{graphicx}

\usepackage[numbers,sort&compress]{natbib}

\begin{document}

%\title[AI-assisted Math Grading]{Exploring AI-assisted Automated Short Answer Grading of Handwritten Mathematics Exams at University Level}
\title[AI-assisted Math Grading]{AI-assisted Automated Short Answer Grading of Handwritten University Level Mathematics Exams}

\author[1]{\fnm{Tianyi} \sur{Liu}
\email{tianliu@student.ethz.ch}}

\author[2]{\fnm{Julia} \sur{Chatain}~\orcidlink{0000-0003-1626-4601}}
\email{julia.chatain@sec.ethz.ch}

\author*[3]{\fnm{Laura} \sur{Kobel-Keller}~\orcidlink{0000-0001-6857-9870}}
\email{laura.kobel-keller@math.ethz.ch}

\author[4,5]{\fnm{Gerd} \sur{Kortemeyer}~\orcidlink{0000-0001-6643-9428}}
 \email{kgerd@ethz.ch}

\author[3]{\fnm{Thomas} \sur{Willwacher}}
\email{thomas.willwacher@math.ethz.ch}

\author[1]{\fnm{Mrinmaya} \sur{Sachan}~\orcidlink{0000-0001-8787-8681}}
\email{mrinmaya.sachan@inf.ethz.ch}

\affil[1]{\orgdiv{Department of Computer Science}, \orgname{ETH Zurich}, \orgaddress{\street{Rämistrasse 101}, \postcode{8092} \city{Zurich}, \country{Switzerland}}}
\affil[2]{\orgname{Singapore-ETH Centre}, \orgaddress{\street{1 College Ave E, CREATE Tower}, \city{Singapore} \postcode{138602}, \country{Singapore}}}
\affil[3]{\orgdiv{Department of Mathematics}, \orgname{ETH Zurich}, \orgaddress{\street{Rämistrasse 101}, \postcode{8092} \city{Zurich}, \country{Switzerland}}}
\affil[4]{\orgdiv{Rectorate and ETH AI Center}, \orgname{ETH Zurich}, \orgaddress{\street{Rämistrasse 101}, \postcode{8092} \city{Zurich}, \country{Switzerland}}}
\affil[5]{\orgname{Michigan State University}, \orgaddress{\street{426 Auditorium Road}, \city{East Lansing}, \state{Michigan} \postcode{48823}, \country{USA}}}

\abstract{
Effective and timely feedback in educational assessments is essential but labor-intensive, especially for complex tasks. Recent developments in automated feedback systems, ranging from deterministic response grading to the evaluation of semi-open and open-ended essays, have been facilitated by advances in machine learning. The emergence of pre-trained Large Language Models, such as GPT-4, offers promising new opportunities for efficiently processing diverse response types with minimal customization. This study evaluates the effectiveness of a pre-trained GPT-4 model in grading semi-open handwritten responses in a university-level mathematics exam. Our findings indicate that GPT-4 provides surprisingly reliable and cost-effective initial grading, subject to subsequent human verification. Future research should focus on refining grading rules and enhancing the extraction of handwritten responses to further leverage these technologies.
}

\date{\today}% It is always \today, today,

\maketitle

\section{Introduction}

In the realm of university-level mathematics education, providing effective feedback through assessment is crucial~\cite{bransford2000people,houston2001assessing,iannone2013students}. However, due to the scarcity of grading personnel, instructors in high-enrollment courses often have to resort to automatically graded assessments, such as multiple-choice or numerical input questions~\cite{weigand2024mathematics}. This reliance on closed-format solutions typically excludes more nuanced forms of assessment like short-answers, ranging from a few lines of derivations to more complex mathematical reasoning, yet not reaching the complexity of open-ended mathematical proofs. These semi-open assessments offer deeper insights into the students' thinking and reasoning than closed-format answers; however, they come with two major challenges:
\begin{itemize}
    \item Answers can be expressed in mathematically equivalent ways. For example, when asking about the roots of a polynomial, the answers $x_1=1,\ x_2=2i,\ x_3=-2i$; $x\in\{1,2i,-2i\}$; $x_{1,2}=\pm2i,\ x_3=1$; $x_1=1,\ x_2=-\sqrt{-4},\ x_3=\sqrt{-4}$ would be equivalent (setting aside considerations of style). 
    \item The solutions are often submitted in handwritten format to avoid the cumbersomeness of formatting or special input requirements, and currently they also need to be graded by hand, which limits their applicability.
\end{itemize}

Computer-based approaches such as STACK~\cite{stack}
that has MAXIMA~\cite{maxima}
in its background, M{\"o}bius~\cite{moebius} (formerly MapleTA~\cite{mapleta}), or exercises in the spirit of Khan academy~\cite{khan}, can tackle some of the above mentioned issues but require more preparation, special design and programming of assessments.

In this light, Artificial Intelligence (AI) offers a promising solution to these grading challenges, both due to its ability to understand semantic relationships and its ability to perform contextual Optical Character Recognition (OCR). The integration of AI in higher education assessment has been gaining traction, recognized for its efficacy in handling diverse and semantically varied responses~\cite{crompton2023,zhang2023}, including mathematics education~\cite{gandolfi2024gpt,zhang2024math}. Particularly, Automated Short Answer Grading (ASAG) technologies hold the promise of enhancing assessment capabilities in educational settings with large student populations, including university mathematics~\cite{burrows2015,haller2022survey,michael2024automatic}.

The development of ASAG methods has shown remarkable progress. Early works up to 2015, reviewed by Burrows, Gurevych, and Stein, focused on hand-engineered features~\cite{burrows2015}. More recent advancements, as surveyed by Haller et al. up to 2022, demonstrate a shift towards models that learn representations from large text corpora~\cite{haller2022survey}. Modern ASAG approaches, such as those employed by Smalenberger et al., specifically tailor models to grading tasks, necessitating extensive fine-tuning~\cite{michael2024automatic}. In contrast, recent general-purpose Large Language Models (LLMs) like GPT-4 and Gemini are designed to adapt to a variety of tasks with minimal customization~\cite{gpt4,gemini}. These models, when evaluated using standard ASAG datasets, show performance levels that rival or exceed those of specialized models from just five years prior to the publication of GPT-4~\cite{kortemeyer2024performance}. They have also started to prove their viability in grading Science and History short answers quizzes without prior training~\cite{henkel2024can}. This comparison raises important considerations regarding the efficiency and practicality of general-purpose AI tools for grading intricate short-answer mathematical responses versus more traditional, custom-designed models.

In the context of practical applicability of such tools at technical universities, assessment of STEM topics is most relevant. Our previous work focused on physics~\cite{kortemeyer2023toward,kortemeyer2024grading}, while our current study addresses another fundamental subject area: introductory mathematics. Virtually all students at technical universities need to pass courses on calculus in their first year. The goal is not necessarily a system that completely takes over grading, but a system that assists humans; even if it takes over a sizeable chunk of the workload while correctly identifying those student solutions that need to be evaluated by a human, a lot is gained. Our results in both physics and mathematics are encouraging that LLMs may play an important role in providing frequent and meaningful assessment for our learners, yet also raise concerns about trust. Future work will need to focus on establishing reliable confidence measures to enable a productive collaboration between AI and human graders.

\section{Background}
\subsection{Optical Character Recognition}
Recognizing Handwritten Mathematical Expressions (HME) remains a challenging task owing to the ambiguities
in handwriting input and a strong dependency on contextual information~\cite{hamad2016detailed,iskandar2023application}. Most research in this area takes
a sequence-to-sequence approach, with two popular base methods being Watch, Attend and Parse (WAP)~\cite{zhang2017watch}
and Track, Attend and Parse (TAP)~\cite{zhang2018track}. The former employed a fully-convolutional network as the watcher
to encode handwritten images and a recurrent neural network as the parser to generate a LaTeX sequence,
while the latter used a stack of Gated Recurrent Units (GRUs) as the tracker to model the input strokes and
a GRU with Guided Hybrid Attention (GHA) as the parser. The two base methods performed particularly
poorly on highly complex HMEs consisting of a large number of strokes, prompting the development of
new methods. There is also closed-source commercial software for HME recognition, among
which Mathpix~\cite{mathpix} is one of the most well-known. Mathpix specializes in scanning images of HMEs and
converting them into industry standards like LaTeX, supporting equations, diagrams and tables. Another
model worth mentioning is GPT-4 with vision (GPT-4V)~\cite{gpt4v}, which enables users to instruct GPT-4 to
analyze image inputs. Although it was not specifically trained for recognizing HMEs, it can be used for this
specific use case.
\subsection{Automated Grading}
The initial idea of utilizing machines for educational assessment dates back to the late 1920s when Sidney Pressey introduced the ``Automatic Teacher'' a device capable of posing and grading multiple-choice questions autonomously~\cite{petrina2004sidney}. This concept evolved significantly by the 1960s, with computers connecting to mainframes via Teletypes to grade numerical responses based on exact matches with predefined correct answers~\cite{suppes1966arithmetic}. The proliferation of the internet further enhanced these systems, enabling the grading of more complex types of responses, such as rankings, mix-and-match problems, and numerical answers that include tolerances and physical units. Additionally, systems could handle algebraic solutions based on symbolic mathematical equivalence and string responses evaluated through exact matching or regular expressions~\cite{sangwin2007assessing,kortemeyer2008experiences,jonz1990another,chapelle1990cloze}. When evaluating mathematical answers, this usually involves Computer Algebra Systems (CASs), such as Maxima~\cite{maxima}, or symbolic mathematics packages such as SymPy~\cite{sympy}, embedded into grading platforms such as STACK~\cite{stack} or LON-CAPA~\cite{loncapa}, to deal with mathematical equivalency or evaluate an answer against criteria such as having particular mathematical properties. These deterministic algorithms can reliably grade with high accuracy, making them suitable for high-stakes summative assessments. However, despite their reliability in assessing responses against set criteria, their psychometric reliability in measuring actual learning outcomes remained questionable~\cite{lord2008statistical,brown2013my,kortemeyer2014extending}. Also, this class of systems is not suitable for questions requiring
a combination of mathematical expressions and natural languages. Another drawback is that students generally
do not have much freedom when entering their answers on a computer because the input has to follow a
certain format. In any case, formulas are inherently two-dimensional constructs, and translating them into a
computer-readable string of alphanumeric characters introduces additional cognitive load in already stressful
exam situations that is unrelated to the learning objectives of most introductory STEM courses.

In contrast, completely open-ended questions promote the unrestricted expression of ideas and methods, not confining students to a specific response structure~\cite{pate2012open}. These questions often require essays or extensive mathematical derivations, and the former have been increasingly graded by artificial intelligence systems~\cite{ruseti2024automated}. For the latter, some preliminary work was done with synthetic solution derivations of physics problems~\cite{kortemeyer2023toward} as well as real-world physics exams~\cite{kortemeyer2024grading}. Although these systems provide valuable formative feedback, their accuracy and reliability are generally lower, which diminishes their suitability for high-stakes summative assessments~\cite{kortemeyer2023toward,jamil2023toward,jackson2023trust,conijn2023effects}.

A compromise between the extremes of closed and open-ended questions is provided by semi-open-ended questions, commonly assessed through Automated Short Answer Grading (ASAG) systems~\cite{zhang2022automatic}. These systems, which generally expect one to three sentence responses, must discern paraphrasing and equivalent meanings—a task requiring advanced machine learning techniques~\cite{leacock2003c}. Unlike essays or mathematical proofs, these responses do not encompass an extensive range of potential answers, nor do they involve complex, multi-step reasoning or stylistic evaluations. 

To enhance their accuracy, ASAG systems often undergo subject-matter specific training or fine-tuning, and the use of transformer-based models can benefit from generating various reference answers before deployment~\cite{ahmed2022deep,akila2023novel}. 
In one of the earliest studies by Lan et al.~\cite{lan2015mathematical,lan2019mathematical}, mathematical
expressions within student solutions were transformed using an extension of the classic bag-of-words model
and student solutions were graded based on the cluster they belonged to. Still, non-mathematical texts were
omitted in the grading process. After a few years, Erickson et al.~\cite{erickson2020automated} managed to take non-mathematical
texts into account by using a pre-trained Stanford Tokenizer with Global Vectors for tokenization. They then
used traditional machine learning techniques like Random Forest, XGBoost and deep learning algorithms
like LSTM for grading. Cahill et al.~\cite{cahill2020context} countered the issue from another angle by integrating mathematical
features extracted from m-rater (an automated scoring engine developed by Educational Testing Service) with
textual context including syntactic relationships for grading, which was relatively effective for high school
level algebra problems.

Recently, attention shifted to using Large Language Models ``out-of-the-box,'' relying on their pre-training instead of training or tuning them for particular exam topics and grading tasks~\cite{kortemeyer2024performance,schneider2023towards,fagbohun2024beyond}. It is generally found that today's LLMs achieve accuracy levels comparable to specially-trained models just a few years ago, but cannot be trusted for standalone grading without human intervention. Instead, they can be consulted for a first round of grading or a ``second opinion,'' and all the more important, it will be to establish reliable measures of confidence in their grading results~\cite{kortemeyer2024grading}.
In this exploratory study, we regard mathematical expressions as a language and aim to employ a pre-trained Large Language Model ``out-of-the-box'' to grade short answers in university-level mathematics.

\section{Methodology and Early Findings}
We aimed to explore a diverse range of mechanisms and workflows for the automatic grading of handwritten mathematics exams, using different OCR tools, methods of answer extraction, prompting techniques, grading-rule formulations, output formats, sampling mechanisms, and confidence measures. However, consequently evaluating all combinations would have resulted in about 100~scenarios; due to the associated high computational load, we abandoned approaches that showed little promise in initial experiments before grading the full exam. In this section, we describe the initial exploratory process and how we arrived at four combinations for the final comparison in Section~\ref{sec:result}. 

\subsection{Data Collection}
\begin{figure}
\begin{center}
\includegraphics[width=\columnwidth]{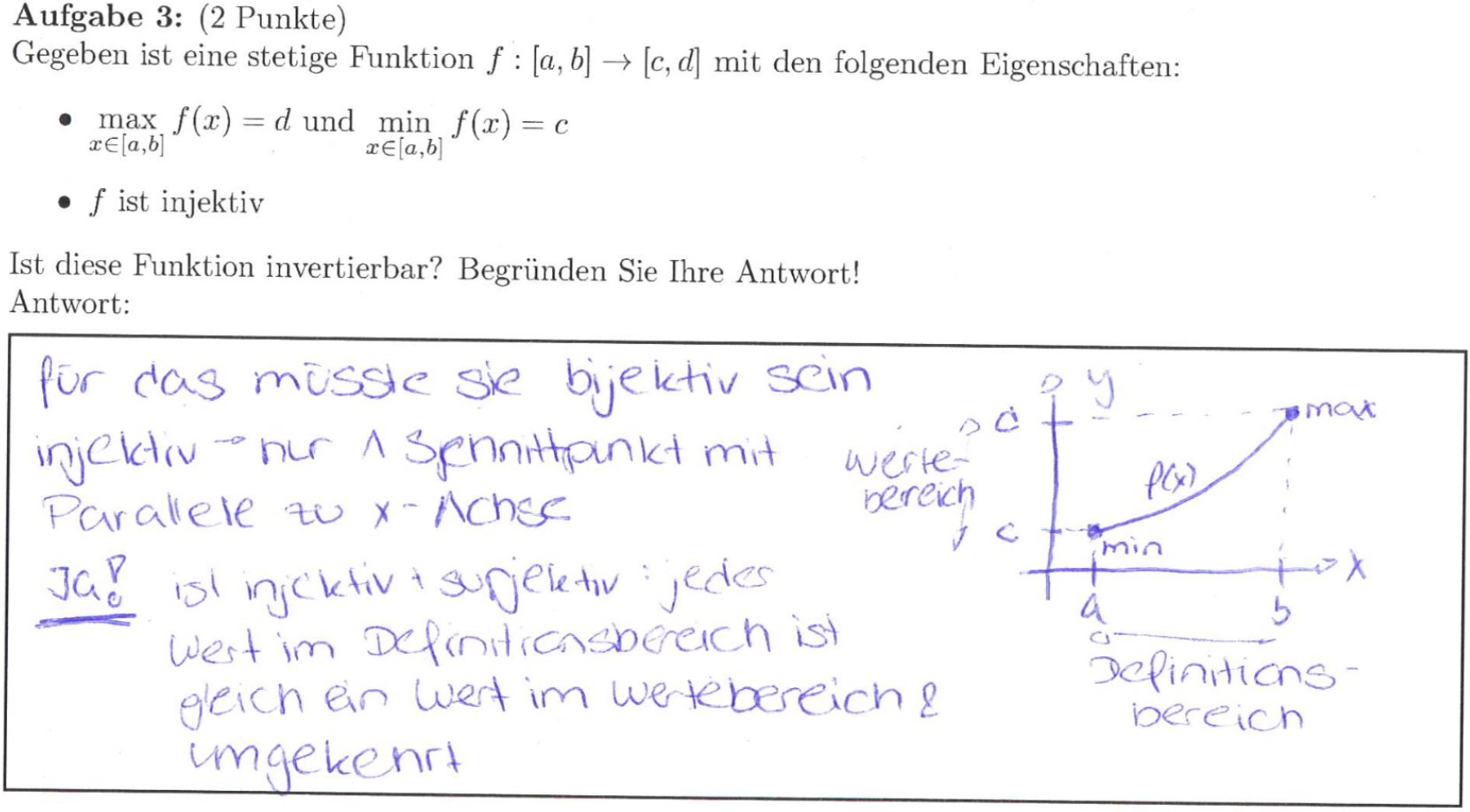}
\end{center}
\caption{Example of an exam problem. Instead of a formula for a function, the student provided a graph.}
\label{fig:ex3}
\end{figure}

\begin{figure*}
\begin{center}
\includegraphics[width=\textwidth]{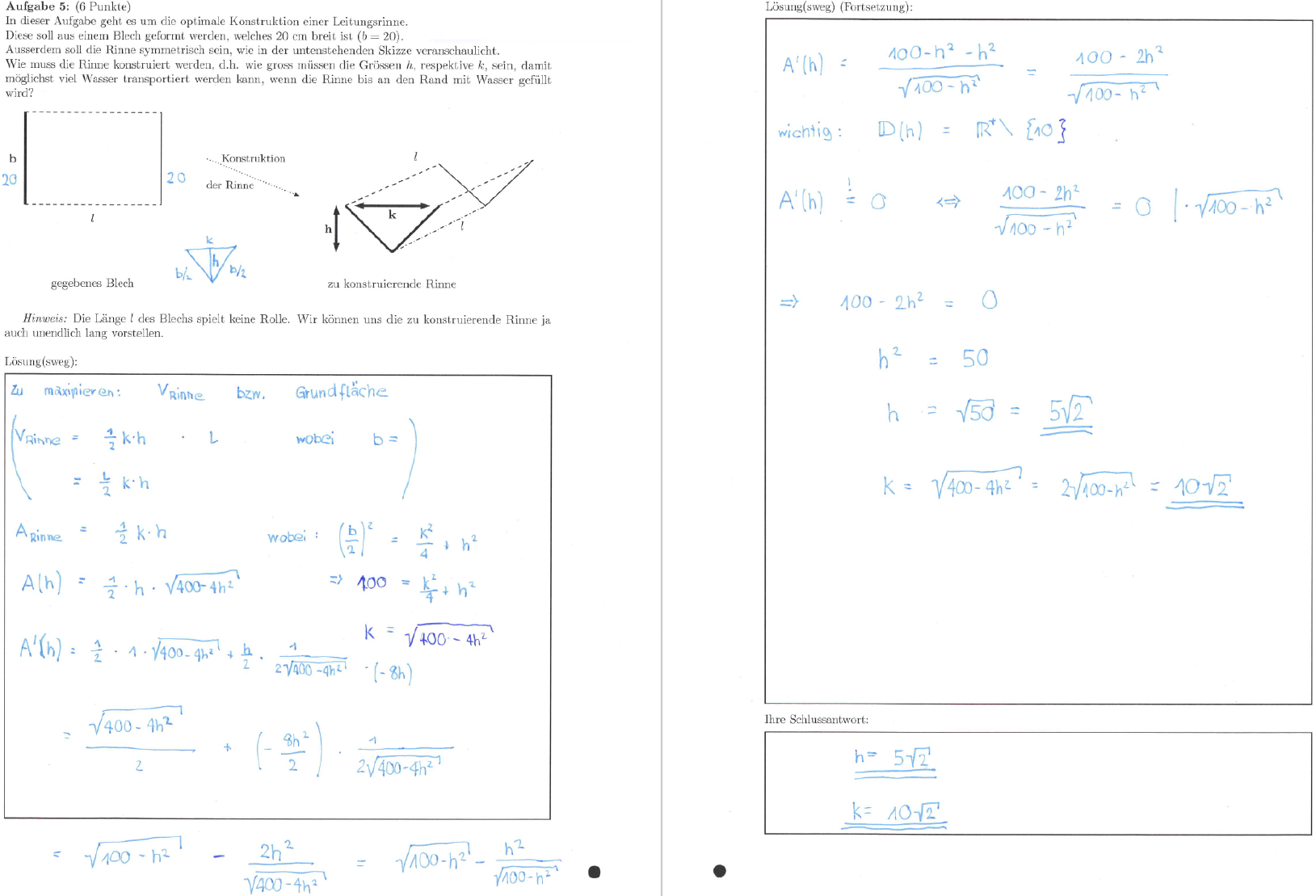}
\end{center}
\caption{Example of an exam problem. The answer extends over several boxes, and the student also wrote outside the box. The last box is the final answer.}
\label{fig:ex5}
\end{figure*}

Our study, approved by the ETH Ethics Commission (proposal EK 2023-N-169), considered data from a mock mathematical exam in German (in this paper, quotations from the exam and its grading were translated to English). This exam consisted of six first-year undergraduate-level questions, totaling nine sub-questions, covering topics such as integrals, Taylor series, optimization, and Ordinary Differential Equations (ODEs). The questions varied in formats, ranging from short numerical answers to composite mathematical reasoning, and integrating various levels of mathematical and natural language. The expected answers thus range from simple numerical responses to detailed justifications and proofs, ranging from one-line mathematical expressions to a handful of lines of symbolic derivations and short-sentence explanations. Figures~\ref{fig:ex3} and~\ref{fig:ex5} shows two of the problems with students answers; the English translation of the whole exam can be found in Appendix~\ref{sec:A1}

\begin{figure}
\begin{center}
\includegraphics[width=\columnwidth]{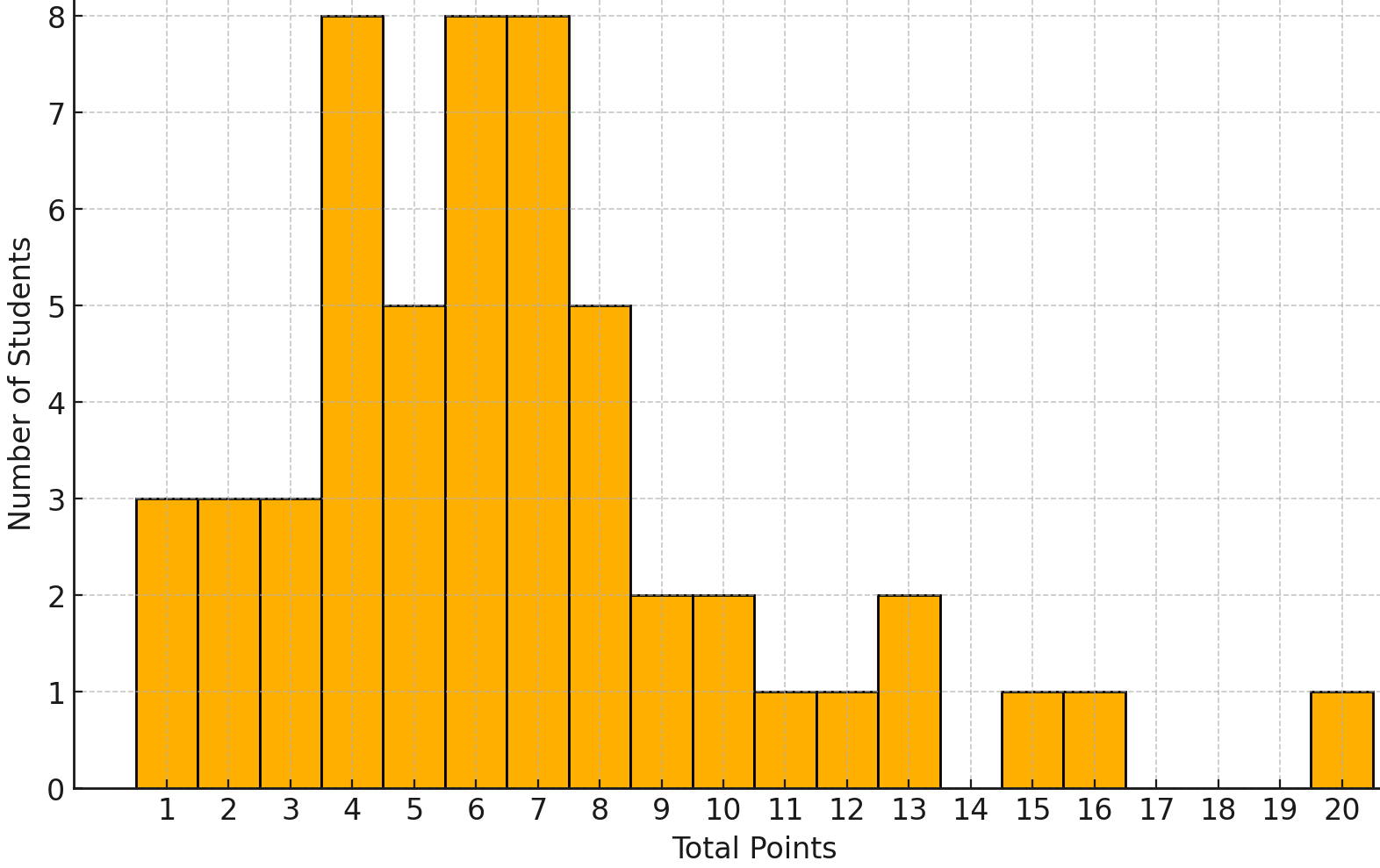}
\end{center}
\caption{Histogram of points (out of 21) achieved by the study participants based on grading by humans (ground truth).}
\label{fig:histo}
\end{figure}

105 students participated in the mock exam, among which 54 consented to the use of their answers. Students could achieve up to 21~points; as Fig.~\ref{fig:histo} shows, the majority of students scored between 4~and 8~points. Notably, there are a few outliers with exceptionally high scores of 15, 16, and 20~points

During the exam, students were instructed to write their answers in designated answer boxes on the answer sheets.
Answer boxes that were detected as ``empty'' during the subsequent workflow were assumed to be unanswered questions and discarded from the analysis, as the assigned points would trivially be zero.

\subsection{Evaluation Methods}
The ground truth for the study was established by human grading. These were also the points that in the end were reported to the students.

Agreement with ground truth was measured in two ways:
\begin{description}
    \item[Accuracy (Acc):] probably the most intuitive measure, describing the observed percentage agreement between the LLM points and the ground truth points.
    \item[Krippendorf's alpha ($\alpha$):] this measure quantifies the observed disagreement corrected for the disagreement expected by chance~\cite{krippendorff2018content}, in our case between the LLM and ground truth. A value of $\alpha\ge0.8$ is usually considered reliable agreement, while $\alpha<2/3$ is usually considered unacceptable. Negative values indicate disagreement that is worse than random chance.
\end{description}

Likely due to the fact that this was only a mock exam, many students might not have sufficiently prepared for it, and the ground truth included a large number of zero-grades. Particularly problems 2,~4, and~6.d had over 90\% solutions that received no points, which corresponds to only four, two, and three students getting credit, respectively.

\begin{figure}
\begin{center}
\fbox{\small
\begin{minipage}{0.96\columnwidth}
{\fontfamily{cmss}\selectfont
{\bf System:}
Extract the text from the image in LaTeX. The output should
only contain the text in LaTeX. If no text is identified in the
image, return Empty.\vspace*{3mm}

{\bf User:}
[Base 64 encoded images]
}
\end{minipage}
}
\end{center}
\caption{Template for the GPT-4V transcription of pre-extracted answer boxes.}
\label{fig:ocrprompt}
\end{figure}

\subsection{Optical Character Recognition}\label{sec:ocrmethods}
After scanning the exams, the second step in the workflow is Optical Character Recognition (OCR). Different approaches were taken to assess their viability:
\begin{description}
    \item[Pre-extracted answer boxes:] Using LaTeX measurements and a Python script, the content of the answer boxes was extracted from the scanned sheets and turned into images. These were interpreted using two different methods:
    \begin{itemize}
        \item In one workflow, we used Mathpix~\cite{mathpix}, a tool to convert mathematical expressions to a variety of formats, including LaTeX.
        \item In a separate workflow, GPT-4V~\cite{gpt4v} was prompted to convert the images to LaTeX; Fig.~\ref{fig:ocrprompt} shows the prompt we used. This was done at two temperature settings, $T=1$ and $T=0$. We also investigated the influence of including the question into the prompt in Fig.~\ref{fig:ocrprompt} as additional context information to transcribe the answer.
    \end{itemize}
    \item[Whole pages:] GPT-4V~\cite{gpt4v} was prompted to convert the complete pages to LaTeX at $T=0.7$, including the questions. The answers were then extracted using GPT-4 and some manual labor (post-extraction).
\end{description}

\subsection{Grading Criteria}
\subsubsection{Grading Rubric}\label{sec:gradrubric}
When teaching assistants grade exams, they are usually given a grading rubric that defines expectations and assignable point values to ensure fairness and consistency. Due to their granularity and particularity, rubrics can contribute to making grading decisions more transparent and reliable~\cite{reddy2010review}.
These rules are written for humans; humans who are knowledgeable about the subject matter in general and the course in particular. The same assumptions cannot be made for a pre-trained general-purpose LLM.

\begin{figure}
\begin{center}
\fbox{\small
\begin{minipage}{0.96\columnwidth}
{\fontfamily{cmss}\selectfont
Grading rule: [Grading rule]

Student answer: [Answer]
}
\end{minipage}
}
\end{center}
\caption{Prompt template for providing grading rules.}
\label{fig:prompt}
\end{figure}

We prompted the model in English, both because GPT-4 tends to perform better in English than in other languages~\cite{achiam2023gpt} and because two of the investigators were not fluent in German. Using one of the workflows as test case, though, it turned out that language had no impact on grading performance beyond random fluctuations; this is maybe not surprising, since the primary language of the documents is ``mathematics.''

In addition to language-dependence, we also tested the rules for robustness by having GPT-4 provide five variations of the grading rules, using the prompt ``\texttt{generate a variation of the following instruction while keeping
the semantic meaning.}''

Fig.~\ref{fig:prompt} shows the template we used for providing the grading rules. Chain-of-thought prompting techniques suggest that demanding an explanation, which would reflect the model thinking process, could improve reasoning and thus grading performance~\cite{wei2022chain}, however, using one workflow as test case, we found no influence on grading beyond random fluctuation. These explanations, however, helped us better understand the grading decisions and possibly improve the rubric items, which is why we kept them.

The italicized statement to ignore irrelevant information was included in the workflows for the pre-extracted answers, but not in the whole page approach. It was primarily designed to address scenarios where additional information was presented by the student that pertains to later grading rules. It also deals with scribbled-out formulas or annotations that student may have made. It is unclear how to handle situations where the student offers additional information that does not appear in any of the grading rules in the case that this information is incorrect; this is also not covered by the grading rules for humans.

We explored two versions of the grading rubric:
\begin{description}
    \item[Original:] Using the same rules that were used by the human graders, even if they included multiple criteria.
    \item[Itemization:] We split grading rules that included multiple criteria for partial credit into multiple finer-grained rules that can be answered with binary judgements. This included scenarios where students could receive full credit for a final correct answer, independent of them listing intermediate steps, or could receive partial credit for intermediate steps even if the final answer is wrong.
\end{description}

\begin{figure}
\begin{center}
\fbox{\small
\begin{minipage}{0.96\columnwidth}
{\fontfamily{cmss}\selectfont
Determine whether the student answer includes the solution
in the grading rule. {\it Ignore the additional information in the student answer that is irrelevant to the grading
rule.}

Provide a short explanation of your
decision.

The output should strictly use the following template:

Judgement: [Yes/No]

Explanation: [Explanation]
}
\end{minipage}
}
\end{center}
\caption{Prompt template for verbalized judgement output with the explored ``ignore''-statement in italics.}
\label{fig:verbalized}
\end{figure}

\begin{figure}
\begin{center}
\fbox{\small
\begin{minipage}{0.96\columnwidth}
{\fontfamily{cmss}\selectfont
Determine whether the student answer includes the solution
in the grading rule. {\it Ignore the additional information in the student answer that is irrelevant to the grading
rule.}

Choose from: (A) Yes (B) No (C) I am not sure.

Provide a short explanation of your decision.

The output should strictly use the following template:

Judgement: [A/B/C]

Explanation: [Explanation]
}
\end{minipage}
}
\end{center}
\caption{Prompt template for multiple-choice output with the explored ``ignore''-statement in italics.}
\label{fig:mcq}
\end{figure}

In addition, the output expected from the LLM may influence its judgement. We considered two different output formats:

\begin{description}
    \item[Verbalized Judgement Format:] This semi-open format expects a verbal output of ``Correct'' or ``Incorrect,'' see Fig.~\ref{fig:verbalized}.
    \item[Multiple Choice Format:] This format is more restrictive, returning the judgement in a multiple-choice format, see Fig.~\ref{fig:mcq}.
\end{description}

Investigating one workflow, we found that in multiple-choice format, the option ``I am not sure'' was never used.
Overall, we found the performance of multiple-choice and verbalized judgements to be very similar, with a slight preference of the verbalized format; we thus discarded the multiple-choice format from later experiments.

\subsubsection{Free Grading}
Another method explored was providing the question, the student answer, and the maximum number of points; no sub-criteria or rubric items regarding different properties of the answer were specified; Fig.~\ref{fig:freeform} shows the associated prompt.

\begin{figure}
\begin{center}
\fbox{\small
\begin{minipage}{0.96\columnwidth}
{\fontfamily{cmss}\selectfont
{\bf System:} 

Based on the question and the maximum number of points,
determine the number of points to be awarded to the student
answer. The number must be an integer. Provide an explanation
for your decision.

The output should strictly use the following template:

Points: [Number of points]

Explanation: [Explanation]
\vspace*{3mm}

{\bf User:}

Question: [Question]

Maximum points: [Number of points]

Student answer: [Answer]
}
\end{minipage}
}
\end{center}
\caption{Prompt template for free grading.}
\label{fig:freeform}
\end{figure}

This approach completely relied on the mathematics knowledge of the LLM, since no master solution was provided. We had found earlier that for some short answer questions, withholding the master solution can improve performance~\cite{kortemeyer2024performance}.

In the context of this mathematics exam, we quickly found that this approach was unpromising, as compared to the thus far best rubric-based grading approach with Acc=$0.69$ and $\alpha=0.60$. Indeed, accuracy for this approach was Acc=$0.50$, $\alpha=0.32$, and Acc$=0.51$, $\alpha=0.34$ when taking a majority or averaging approach to sampling the output, respectively. We thus discarded it in favor of grading rubrics.

\subsection{Sampling}
Given the probabilistic nature of LLMs, our approach leverages the Law of Large Numbers by conducting multiple grading iterations at a temperature $T>0$~\cite{kortemeyer2024grading}. In all scenarios,
we run the same grading prompts five times at
$T = 0.7$. 

For the whole-pages approach (see Sect.~\ref{sec:ocrmethods}), we also run several OCR rounds. The grading process begins with OCR using the Whole Pages approach (see Sect.~\ref{sec:ocrmethods}), which generates several versions of each page. As each version of a student's answer is graded five times, this results in 25~potentially different point values.

For the pre-extracted answer boxes, we
used a majority vote, while for the whole-pages
approach, we used averages.
\subsection{Confidence}\label{sec:confidence}
In realistic grading scenarios, where AI is supposed to replace some of the human grading effort, there would not be a ground truth to fall back on. It is thus very important to decide which AI results can be trusted, and which AI results need additional human attention. We attempt to address this by establishing a confidence measure.

For a probabilistic approach to confidence determination, we prompted the system for numerical responses (``points'') and used averages and standard deviations. We then based the final decision on whether to trust the AI-grading on the mean and standard deviation of these point values. A small standard deviation indicates high confidence and presumably alignment with the ground truth grade, while a large standard deviation would be associated with uncertainty.

The grading would result in the assignable point value nearest to the mean. If at most one assignable point value falls within the standard deviation $\sigma$ of the mean, the grading is categorized as ``Can Decide,'' otherwise, it is labeled as ``Cannot Decide.'' For example, available point values might be $0.0$, $0.5$, $1.0$, $1.5$, and $2.0$; if the system arrives at $1.21\pm0.13$ or $1.21\pm0.23$, a point value of $1.0$ would be assigned with certainty as no or one assignable point value is inside of the standard deviation, respectively, while $1.21\pm0.31$ would be a case of ``Cannot Decide'' as two assignable point values fall within the standard deviation $\sigma$. Based on those, a contingency table can be set up, see Table~\ref{tab:contingency}.

\begin{table}
\centering
\begin{tabular}{r|ll}
                   &can decide         &cannot decide\\\hline
correct    &true positive (TP) &false negative (FN)\\
incorrect &false positive (FP)&true negative (TN)
\end{tabular}

\caption{Contingency table for the probabilistic confidence measure; the correctness of the judgement is based on the ground truth.}

\label{tab:contingency}

\end{table}

All instances labeled as ``Cannot Decide'' would require evaluation by teaching assistants. False negatives, where the system could not decide but the rounded grade matches the correct grade, create extra work but no harm. More critical are false positive cases, where the system incorrectly confirms a grade, as these represent the only detrimental outcomes.

\section{Results}\label{sec:result}

In this section, we describe our main results regarding the three main steps of the grading workflow: answers extraction, prompting, and actual grading.

\subsection{Answer Extraction}
\begin{figure}
\begin{center}
\includegraphics[width=\columnwidth]{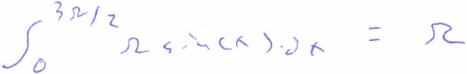}
\end{center}
\caption{Example of a pre-extracted student answer to problem~1 (see App.~\ref{sec:A1}).}
\label{fig:studanswer}
\end{figure}

Fig.~\ref{fig:studanswer} shows an example of an extracted student answer.
Answer extraction proved to be more challenging than initially anticipated, particularly since students tended to write across box boundaries and even using the margins of subsequent problems.

Students left several questions unanswered, but as Table~\ref{tab:empty_answers} shows, extracting the answer boxes first resulted in more unrecognized answers than were actually left empty. In this scenario, Mathpix detects more answers than GPT-4V. For GPT-4V, apparently more information is preserved at $T=1$ than at $T=0$. In particular, we found that GPT-4V is more likely to recognize very short mathematical expressions inside the answer
boxes at $T=1$, but ignores them at $T=0$. As noted earlier, unrecognized answer boxes were assumed to be unanswered questions and discarded from the analysis; with only 54 students in the study, this led to low sample sizes, particularly for question~6.d. 

As higher temperatures lead to better results, we performed the interpretation of whole pages at $T=0.7$. This post-extraction mechanism was generally more reliable, however, for some problems it surprisingly detected more answers than were actually provided by the students; as it turned out, in those cases it interpreted the question text as an answer. While more likely than not, the question text will not be graded as a correct answer later in the workflow, this could likely have been avoided by embedding textual rather than graphical markers into the exam text.

\begin{table*}
\centering
\begin{tabular}{lr|rrrr||r}
Problem&Ground&\multicolumn{3}{c}{Pre-extracted answer boxes}&Whole pages&Ground\\
&Truth&Mathpix&GPT-4V& GPT-4V&GPT-4V&Truth\\
&Empty Box&&$T=1$&$T=0$&$T=0.7$&Zero Grade\\\hline
1 & 6\% & 19\% & 35\% & 67\% & {\bf 6\% }& 65\%\\
2 & 20\% & 31\% & 33\% & 74\% & {\bf 15\% }& 93\%\\
3 & 6\% & 24\% & 17\% & 74\% & {\bf 7\% }& 31\%\\
4 & 35\% & 46\% & 44\% & 61\% & {\bf 24\% }& 95\%\\
5 & 2\% & 11\% & 11\% & 11\% & {\bf 2\% }& 11\%\\
6.a & 2\% & 22\% & 19\% & 48\% & {\bf 6\% }& 6\%\\
6.b & 26\% & 44\% & 41\% & 54\% & {\bf 19\% }& 67\%\\
6.c & 33\% & 46\% & 46\% & 54\% & {\bf 26\% }& 76\%\\
6.d & 61\% & 70\% & 70\% & 89\% & {\bf 57\% }& 94\%\\
\end{tabular}

\caption{Comparison of percentage of empty (ground truth) and unrecognized answer boxes for the workflows described in Sect.~\ref{sec:ocrmethods} (see Appendix~\ref{sec:A1} for the problems themselves). Note that unrecognized answer boxes were discarded from the analysis. Also included in the rightmost column is the percentage of zero-grades in the ground truth. Boldface values indicate the lowest percentage of unrecognized answers for each problem.}
\label{tab:empty_answers}
\end{table*}

For pre-extracted answer boxes, we also explored including the question into the prompt in Fig.~\ref{fig:ocrprompt}. For the example in Fig.~\ref{fig:studanswer}, transcription without the question resulted in
\[\int_0^{3\pi/2}\Omega(x(t))\cdot x(t)\ dt = \Omega\ ,\]
while including the text of problem~1 in App.~\ref{sec:A1} resulted in
\[\int_0^{3\pi/2}\pi\sin(x)\ dt = \pi\ .\]
Even for a human, the student answer would have been hard to decipher without background information, and for the LLM, its inclusion clearly improved the accuracy of the transcription in this case. There is, however, also a strong caveat that kept us from generally applying this mechanism: previous experiments found anecdotal evidence that the same practice can also lead to the LLM ``fixing'' the student answer, i.e., ``seeing what it expects.'' 

\subsection{Robustness of Grading Prompts}
To test robustness of AI grading with respect to wording and phrasing of the grading rules, we used GPT-4 to generate five variations of the itemized grading rules (see Sect.~\ref{sec:gradrubric}). Since the resulting grades from the five variations turned out to have a Krippendorff’s
alpha of 0.88  among each other, GPT-4’s grading appears generally robust against paraphrasing of the grading prompt as long as the meaning of the rule remains unchanged.

\subsection{Probabilistic Confidence Measure}
Table~\ref{tab:prop} shows the result of applying the $\sigma$-confidence measure to the probablistic grading approach described in Sect.~\ref{sec:confidence}. Here, we used dichotomous contingency measures, which we evaluated in terms of standard measures, namely the accuracy reflecting the percentage agreement,
the precision reflecting the accuracy of correctness, the recall reflecting the completeness of correct grades, the harmonic mean F1 of precision and recall, and the percentage of false positive cases.

\begin{table*}
\centering
\begin{tabular}{rrrrrr}
Problem & Accuracy & Precision & Recall & F1 & FP rate \\\hline
1 & 0.83 & 0.80 & 1.00 & 0.89 & 0.17 \\
2 & 0.93 & 0.92 & 1.00 & 0.96 & 0.07 \\
3 & 0.57 & 0.57 & 0.97 & 0.72 & 0.41 \\
4 & 0.94 & 0.94 & 1.00 & 0.97 & 0.06 \\
5 & 0.44 & 0.13 & 0.25 & 0.17 & 0.39 \\
6.a & 0.61 & 0.58 & 1.00 & 0.73 & 0.39 \\
6.b & 0.72 & 0.75 & 0.92 & 0.83 & 0.22 \\
6.c & 0.93 & 0.93 & 1.00 & 0.96 & 0.07 \\
6.d & 0.93 & 0.93 & 1.00 & 0.96 & 0.07 \\
\end{tabular}
\caption{Performance of the confidence measure based on the standard deviation $\sigma$. Note that accuracy here is that of the confidence measure, not that of the grading. The FP-rate describes the percentage of grades deemed reliable by the $\sigma$-confidence measure, but different from the ground truth.}
\label{tab:prop}
\end{table*}
As can be seen, the accuracy of this criterion in correctly identifying AI-grades that correspond to the human-assigned grades varies greatly between $0.44$ and $0.94$. It turns out that high confidence accuracy is, with the exception of item 6.c, associated with problems that have high failure rates (zero grade) in Table~\ref{tab:empty_answers}; removing problems~2,~4 and~6.d from consideration, the average accuracy of the $\sigma$-confidence criterion is $0.69$. More importantly, the average FP-rate is $0.27$; in other words, for the problems were we had sufficient statistics, in about a quarter of the grading instances, the algorithm indicated confidence into a grading result that did not correspond to the ground truth.

\subsection{Grading Accuracy}
Problems~2,~4, and~6.d had insufficient sample sizes of correct answers and were discarded from the final analysis; for the remaining items, Table~\ref{tab:compare_accs} shows the agreement measures for the workflows that were identified as the most promising for all remaining problems.

The top row of the header shows the method of answer extraction, the second row the mode of prompting, i.e., the original rubric versus the itemized version of rubric items with multiple criteria (see Sect.~\ref{sec:gradrubric}), and the third row lists the OCR mechanism (see Sect.~\ref{sec:ocrmethods}). For the $\sigma$-confidence column, only the items that were deemed reliable by the algorithm described in Sect.~\ref{sec:confidence} were included in the analysis; the rightmost columns shows the percentage of responses meeting the confidence condition.

\begin{table*}
\centering
\begin{tabular}{rrr|rr|rr|rr|rrr}
& \multicolumn{6}{c|}{Pre-extracted answer boxes} & \multicolumn{5}{c}{Whole pages}\\
& \multicolumn{2}{c|}{original} & 
\multicolumn{4}{c|}{itemized} & \multicolumn{5}{c}{original}\\
& \multicolumn{2}{c|}{Mathpix} & \multicolumn{2}{c|}{Mathpix} & \multicolumn{2}{c|}{GPT-4V}
& \multicolumn{5}{c}{GPT-4V} \\
&&&&&&&&&\multicolumn{3}{|c}{$\sigma$-confidence}\\
Problem & Acc & $\alpha$ & Acc & $\alpha$ & Acc & $\alpha$ & Acc & $\alpha$& Acc & $\alpha$ & Positive\\\hline
1 & {\bf 0.91} & {\bf 0.85} & 0.89 & 0.81 & 0.75 & 0.46 & 0.67 & 0.52 & 0.79 & 0.52 & 84\%\\
3 & 0.27 & 0.15 & 0.41 & 0.06 & 0.49 & 0.16 & 0.54 & 0.60 & {\bf 0.55} & {\bf 0.66} & 94\%\\
5 & {\bf 0.27} & {\bf 0.06} & 0.15 & -0.34 & 0.17 & -0.15 & 0.21 &-0.12 & 0.09 & -0.63 & 43\% \\
6.a & {\bf 0.74} & {\bf 0.52} & {\bf 0.74} & 0.47 & 0.24 & -0.21 & 0.53 & 0.06 & 0.57 & 0.07 & 92\%\\
6.b & 0.70 & 0.36 & 0.60 & 0.24 & 0.60 & 0.41 & 0.68 & {\bf 0.56} & {\bf 0.71} & 0.44 & 83\%\\
6.c & 0.86 & 0.70 & 0.83 & 0.62 & 0.83 & 0.65 & {\bf 0.90} & {\bf 0.76} & {\bf 0.90} & {\bf 0.76} & 100\%
\\\hline
Average & {\bf 0.62} & {\bf 0.44} & 0.60 & 0.31 & 0.51 & 0.22 & 0.59 & 0.40& 0.60 & 0.30& 83\%
\end{tabular}
\caption{Performance of the workflows under investigation. For the probabilistic workflow that discards grading outcomes based on the $\sigma$-criterion, we also list the percentage of grades deemed reliable.}
\label{tab:compare_accs}
\end{table*}

Itemizing the multi-criterion rules generally did not result in the expected performance increase; it is, however, a helpful technique for particular problems, for example the partial-credit rule for  Problem~3. Also the application of the $\sigma$-confidence criterion did not significantly improve performance. In case of problem~5, which has a low or even negative $\alpha$ across all workflows, the accuracy even decreased; while the algorithm was particularly selective in this case (less than half of the grading instances passed the criterion), it tended to pick false positives.

The average results for the original rubric formulations are comparable, with an accuracy range from $0.59$ to $0.62$. Overall, the performance of the grading mechanisms is not acceptable for high-stakes summative assessments. All the more important is the false-positive rate, where the system identified a judgement as reliable when in fact it was not corresponding to the ground truth; this value varies between 7\% and 41\%, where the high values correspond to the problems with inaccurate grading; in the pre-extracted answer workflow, these same problems are flagged by a low and even negative $\alpha$. 

\section{Lessons}
\subsection{Markers instead of Answer Boxes}
The extraction of answer boxes turned out to miss possible solutions (see Table~\ref{tab:empty_answers}), as students frequently wrote across their borders or even worked ``outside the box.''
Even when not extracting the boxes but performing OCR on the whole page, the workflow was frequently derailed by the answer boxes, as they were essentially extraneous graphical elements.

\begin{figure}
\begin{center}
\fbox{\small
\begin{minipage}{0.96\columnwidth}
{\fontfamily{cmss}\selectfont
\ldots\vspace*{3mm}

Problem 2 (2 points):

For which  $x \in \mathbb{R}$ does the following series converge?

\[
\sum_{k=2}^{\infty} k (k-1) (x-2)^{k-2}
\]
\vspace*{3mm}

Solution of Problem 2:

\begin{eqnarray*}
&\uparrow&\\
&\mbox{\it provide more than sufficient space}&\\
&\downarrow&
\end{eqnarray*}

Problem 3 (2 points):

Given is a continuous function $f: [a,b] \rightarrow [c,d]$ with the following properties:
\begin{itemize}
\item $\displaystyle \max_{x \in [a,b]}f(x) = d$ and $\displaystyle \min_{x \in [a,b]}f(x) = c$\vspace*{1mm}
\item $f$ is injective
\end{itemize}
Is this function invertible? Justify your answer!
\vspace*{3mm}

Solution Problem 3:

\ldots
}
\end{minipage}
}
\end{center}
\caption{Example of an exam layout that facilitates OCR.}
\label{fig:layout}
\end{figure}

Overall, it seems preferable to perform OCR on the complete page, including the questions, which can easily be filtered out if textual markers are in place. This would have worked with both Mathpix and GPT-4V.

In retrospect, we should have provided a layout as shown in Fig.~\ref{fig:layout}. The markers ``Problem $N$ ($n$ points):'' and ``Solution of Problem $N$'' will get read by the OCR and can be used to dissect the textual exam files in a straightforward way, for example by using regular expressions. The answer spaces need to be sufficiently large, so students do not start writing on the margins of the next problem.   

\subsection{AI versus Human Grading Rules}
When defining grading rules, it is important to keep in mind that current pre-trained Large Language Models have neither common sense nor the ability to reliably carry out advanced calculations~\cite{gandolfi2024gpt}. 

For mathematical tasks, at least in our experiments with German and English, we found language to have no significant influence on grading accuracy; in separate experiments, we found that explicit translation steps are not necessary, as GPT-4 fluidly moves between the languages~\cite{kortemeyer2024grading}; unfortunately, this may not apply to all languages, as performance of GPT-4 rapidly falls off with less commonly spoken languages~\cite{achiam2023gpt}. 

Oftentimes, providing the question text alongside the grading rules can provide additional leeway even if the grading rules are imprecise or incomplete, but in any case, grading rules need to be more carefully crafted than what is typically provided for teaching assistants; Figure~\ref{fig:polynomial} shows a typical example for what can happen when giving a grading rule meant for humans directly to GPT-4.

\begin{figure}
\begin{center}
\fbox{\small
\begin{minipage}{0.96\columnwidth}
{\fontfamily{cmss}\selectfont
Provided in System Prompt
\begin{description}
\item[Question:] Determine all roots of the following polynomial:
\[p(x)=x^3-x^2+4x-4\]
\item[Grading rule:] $x = 1$ is a zero of the given polynomial.
\item[Student answer:] $x_0=1,\ 4i,\ -4i$
\end{description}

System Response
\begin{description}
\item[Judgement:] Incorrect
\item[Explanation:] The student has correctly identified $x=1$ as a
root of the polynomial but the other two roots given by the
student are incorrect. The correct roots of the polynomial
$x^3-x^2+4x-4$ are $x=1$, $x=2-\sqrt{3}i$ and $x=2+\sqrt{3}i$.
\end{description}
}
\end{minipage}
}
\end{center}
\caption{Example of providing the questions and imprecise grading rules.}
\label{fig:polynomial}
\end{figure}

The grading was supposed to give one point for having found the real answer; the imaginary answers would be dealt with in a later grading rule. Due to the formulation of the grading rule, GPT-4 attempted to determine the additional roots by itself, which at the time was a task beyond its capabilities. Here, changing the grading rule to ``judge as correct if the student found $x=1$ as the only real solution; disregard any complex solutions for now'' leads to the desired result: ``Given that the grading rule focuses only on the identification of the real root $x=1$ and explicitly disregards the assessment of complex solutions, the student's answer should be judged as correct concerning the grading rule.'' The following grading rule should mirror this, ``judge as correct if the student found $x=2i$ and $x=-2i$ as the only complex solutions; disregard any real solutions.''

\section{Limitations}
This study clearly has some limitations. The exam for data
collection was not mandatory and students might not have exerted their best effort, resulting in a high
percentage of zero grades. The non-compulsory nature of the exam might also influence students’ handwriting,
potentially making it more difficult for OCR tools to interpret. 

Moreover, due to the language
proficiency of those involved in the project, the grading rules were translated from German to English while
the student answers generally remained in German. Although the language seemed to have little impact on
the experiments, keeping the consistency of the languages would enhance the reliability of the results.

Due to the high computational effort involved, the effect to the temperature could not be fully explored.

\section{Discussion}
The AI-grading in this study was backed up by a ground truth, which in a real-world scenario would defeat the purpose. Instead, the system would need to reliably establish confidence measures, so only those student solutions where the system confidence is low get graded by teaching assistants.
While certainly more work should be done on handwriting recognition and grading rules, a high false-positive rate also indicates the need for developing more reliable confidence measures.

\section{Outlook}
Our preliminary results indicate that GPT-4o~\cite{gpt4o} is superior to GPT-4V in interpreting handwritten mathematical expressions~\cite{kortemeyer2024grading}. Future studies should involve a layout as shown in Fig.~\ref{fig:layout} and interpreting complete pages with GPT-4o before separating out the answers. Finally, considering mandatory exams would likely provide more realistic data.

Confidence measures are essential for building trustworthy grading applications. The application of psychometric methods such as Classical Test Theory or Item Response Theory could be used in identifying outliers and iterative improvement of grading rules~\cite{lord2008statistical,kortemeyer2019quick}.
    
\section{Conclusions}
This study evaluates the performance of GPT-4 in grading undergraduate-level mathematical tasks and
proposes a method to estimate the reliability of GPT-4 generated grades. Based on the observations of
GPT-4’s grading behaviour and issues discovered during the exploration process, recommendations for
formulating grading rules for AI-assisted grading are also offered. Despite the fact that the grading rules used
in the experiments are not curated for this specific context, GPT-4 showcases considerable potential to serve
as an assistant grader to human graders with its natural language understanding and mathematical reasoning
capabilities, especially for short response questions. Its grading is also generally robust against paraphrasing
of the grading prompt when the sentence structure remains consistent. However, its ability in handwritten
mathematical expression recognition is limited, which indicates the need for further effort to improve the
extraction of student answers. Here, we found it most reliable to first scan and transcribe the whole exam sheet and extracting the answers later, in the future possibly delineated by simple text markers.

As to confidence estimation, the novel approach grounded in probability
theories yields promising results, opening up new avenues for further research in this direction. Although not
covered in this work, fine-tuning-based calibration is another approach worth considering. However, additional work on confidence measures is crucial for making AI-assisted grading for high-stakes mathematics exams a viable option.
\section*{Declarations}

\subsection*{Funding}
No funding was received.
\subsection*{Conflict of interest/Competing interests}
No authors have a conflict of interest or competing interests.
\subsection*{Ethics approval and consent to participate}
The study was approved by the ETH Ethics Commission as proposal EK 2023-N-169.

\subsection*{Consent for publication}
All authors gave their consent to publication. 

\subsection*{Data availability} 
Upon request, subject to an amendment of the human-subject research protocol.

\subsection*{Materials availability}
Upon request.

\subsection*{Code availability}
Upon request.

\subsection{Author contribution}
LKK and TW conceptualized the study. LKK wrote and graded the exam. TW provided the LaTeX and Python code for answer-field extraction. TL wrote the code for the workflow and carried out the data analysis with JC and GK as advisors and MS as academic supervisor. JC provided the technical expertise for the usage of the Mathpix and GPT-4(V) APIs. GK wrote the initial version of the manuscript using materials from TL's thesis, and all authors contributed to the submitted manuscript.

\begin{appendices}

\section{English Translation of the Exam}\label{sec:A1}

{\bf \noindent Problem 1:} (2 points)\\
\begin{displaymath}
\int_{0}^{3\pi/2}\pi \sin (x) \; dx = \; ?
\end{displaymath}

{\bf \noindent Problem 2:} (2 points) \\
For which $x \in \mathbb{R}$ does the following series converge?
\begin{displaymath}
\sum_{k=2}^{\infty} k (k-1) (x-2)^{k-2}
\end{displaymath}

{\bf  \noindent Problem 3:} (2 points) \\
Given is a continuous function $f: [a,b] \rightarrow [c,d]$ with the following properties:
\begin{itemize}
\item $\displaystyle \max_{x \in [a,b]}f(x) = d$ and $\displaystyle\min_{x \in [a,b]}f(x) = c $
\item
$f$ is injective
\end{itemize}
Is this function invertible? Justify your answer!\\

{\bf  \noindent Problem 4:} (2 points) \\
Consider the function $f(x)= \sin(\pi x) + 3x^2$. What is the second degree Taylor polynomial when developed around the point $x=1$? \\

{\bf  \noindent Problem 5:} (6 points) \\
This problem involves the optimal construction of a drainage channel. \\
It should be formed from a sheet of metal that is $20$ cm wide ($b=20$). \\
Furthermore, the channel should be symmetric, as illustrated in Fig.~\ref{fig:drain}. \\
What dimensions must the channel have, i.e., what should be the sizes of $h$, and $k$, respectively, to maximize water transport when the channel is filled to the brim?
\begin{figure*}
\begin{tikzpicture}[scale=1.35]
\tikzset{>=latex}
%\draw (0,0) rectangle (5,2);
\draw[ultra thick] (2,0) -- (02,2);
\draw (1.9,1) node[left]{$\textbf{b}$};
\draw (3.5,-0.1) node[below]{$l$};
\draw[dashed] (2,0) -- (5,0);
\draw[dashed] (2,2) -- (5,2);
\draw[] (5,0) -- (5,2);
\draw (3.5,-1.5) node[below]{$\text{provided sheet metal}$};

\draw[dotted,->] (6,1) -- (8,0.2);
\draw (7,0.85) node[above]{$\text{\small{construction}}$};
\draw (6.5,0.27) node[above]{$\text{\small{of drain}}$};

\draw[ultra thick] (9,0) -- (10,-1);
\draw[ultra thick] (10,-1) -- (11,0);

\draw[ultra thick, <->] (8.8,-1) -- (8.8,0);
\draw (8.8,-0.5) node[left]{$\textbf{h}$};

\draw[ultra thick, <->] (9.2,0) -- (10.8,0);
\draw (10,0) node[below]{$\textbf{k}$};

\draw[dashed] (11,0) -- (13,1);
\draw[dashdotted] (10,-1) -- (12,0);
\draw[dashed] (9,0) -- (11,0.9);
\draw[] (11,0.9)  -- (12,0);
\draw[] (13,1)  -- (12,0);
\draw (10.5,1.25) node[below]{$l$};
\draw (11.5,-0.25) node[below]{$l$};

\draw (10,-1.5) node[below]{$\text{to be constructed drain}$};
\end{tikzpicture}
\caption{Exam illustration of the drain pipe assembly.}
\label{fig:drain}
\end{figure*}

\textit{Note:} The length $l$ of the sheet does not matter. We can imagine the channel being infinitely long. \\
\ \\
  
{\bf  \noindent Problem 6:} (7 points) \\
This problem concerns complex numbers and differential equations
\begin{itemize}
\item[a)] (3 points) \\
Determine all roots of the following polynomial
\begin{displaymath}
p(x) = x^3 -x^2 +4x-4.
\end{displaymath}

\item[b)] (2 points) \\
Now consider the following differential equation
\begin{eqnarray*}
&&y^{'''}(x) -y^{''}(x) +4y'(x)-4y(x)\\
&=&y^{(3)}(x) -y^{''}(x) +4y'(x)-4y(x)=0
\end{eqnarray*}
Assuming a solution of the form $y(x)=e^{\lambda x}$, what must hold for $\lambda$? Describe, simplify and justify your answer as precisely as possible! \\

\item[c)] (1 point) \\
Clearly, $y(x)= e^x$ is a solution to the differential equation given in part b). \\
Is the function $f_2(x)= \cos (2x)$ a solution to the differential equations given in part b)? Justify your answer! \\

\item[d)] (1 point) \\
What other function could be a solution to the given differential equation? \\
Consider the analogy between the given differential equation and the polynomial in part a)! \\
\end{itemize}

\end{appendices}

\bibliographystyle{unsrt}
\bibliography{mathgrading}

\end{document}